\documentclass[12pt,leqno,draft]{article}

\newtheorem{theorem}{Theorem}
\newtheorem{lemma}[theorem]{Lemma}
\newtheorem{proposition}[theorem]{Proposition}
\newtheorem{definition}[theorem]{Definition}
\newtheorem{corollary}[theorem]{Corollary}

\newcommand{\begintheorem}{\addtocounter{equation}{1}\begin{theorem}}
\newcommand{\beginlemma}{\addtocounter{equation}{1}\begin{lemma}}
\newcommand{\beginproposition}{\addtocounter{equation}{1}\begin{proposition}}
\newcommand{\begindefinition}{\addtocounter{equation}{1}\begin{definition}}
\newcommand{\begincorollary}{\addtocounter{equation}{1}\begin{corollary}}

\begin{document}

\title{Some notes about matrices, 4}

\author{Stephen William Semmes 	\\
	Rice University		\\
	Houston, Texas}

\date{}

\maketitle

	As usual, we let ${\bf R}$ denote the real numbers, ${\bf C}$
the complex numbers, ${\bf R}^m$ the space of $m$-tuples of real
numbers, and ${\bf C}^n$ the space of $n$-tuples of real numbers.
Also, ${\bf Z}$ denotes the integers, and ${\bf Z}^m$ the standard
integer lattice in ${\bf R}^m$.  We write ${\bf Z}[i]$ for the
\emph{Gaussian integers}, which are the complex numbers of the form
$a + i \, b$ with $a, b \in {\bf Z}$, and $({\bf Z}[i])^n$ for
the standard integer lattice in ${\bf C}^n$.

	It will be convenient to write $\mathcal{L}({\bf R}^m, {\bf
C}^n)$ for the space of real-linear mappings from ${\bf R}^m$ to ${\bf
C}^n$.  The complex structure on ${\bf C}^n$ is still relevant for
this space, in that $\mathcal{L}({\bf R}^m, {\bf C}^n)$ is naturally a
complex vector space.  This is because one can multiply elements of
$\mathcal{L}({\bf R}^m, {\bf C}^n)$ by $i$, and these linear
transformations can be described by $m \times n$ matrices of complex
numbers in the usual manner, using the standard bases for ${\bf R}^m$
and ${\bf C}^n$.

	Let us write $\mathcal{L}^*({\bf R}^m, {\bf C}^n)$ for the
subset of $\mathcal{L}({\bf R}^m, {\bf C}^n)$ consisting of linear
transformations whose kernels are trivial, at least when $m \le 2n$,
so that this is possible.  Using the usual Euclidean topology for
$\mathcal{L}({\bf R}^m, {\bf C}^n)$, $\mathcal{L}^*({\bf R}^m, {\bf
C}^n)$ is an open set.  When $m = 2n$, $\mathcal{L}^*({\bf R}^m, {\bf
C}^n)$ consists of the invertible real-linear transformations
from ${\bf R}^m$ onto ${\bf C}^n$.

	We can define a \emph{lattice} in ${\bf C}^n$ to be the image
of ${\bf Z}^{2n}$ under an element of $\mathcal{L}^*({\bf R}^{2n},
{\bf C}^n)$.  For such a lattice $L$, we get a quotient ${\bf C}^n /
L$ in the usual manner, which is both an abelian group under addition
and a complex manifold, in fact a complex affine manifold.  In other
words, there are nice complex affine local coordinates for the
quotient space coming from those of ${\bf C}^n$.

	Of course the standard lattice $({\bf Z}[i])^n$ is a
lattice in ${\bf C}^n$ in this sense.  For the moment let us restrict
our attention to lattices $L$ in ${\bf C}^n$ which are of the form
$A(({\bf Z}[i])^n)$ for some invertible complex-linear mapping
$A$ on ${\bf C}^n$.  This is a special case, in the same way that
invertible complex-linear mappings on ${\bf C}^n$ are a special
case of elements of $\mathcal{L}^*({\bf R}^{2n}, {\bf C}^n)$.

	Let us write $GL({\bf C}^n)$ for the group of invertible
complex-linear transformations on ${\bf C}^n$, and $SL({\bf C}^n)$ for
the subgroup of $GL({\bf C}^n)$ consisting of linear transformations
with determinant equal to $1$.  Also, we write $U({\bf C}^n)$ for the
group of \emph{unitary} linear transformations on ${\bf C}^n$,
which are the invertible linear transformations which preserve
the standard Hermitian inner product on ${\bf C}^n$, which is the
same as saying that the inverse of the linear transformation is equal
to its adjoint.  The subgroup $SU({\bf C}^n)$ of $U({\bf C}^n)$
consists of the unitary linear transformations on ${\bf C}^n$
which also have determinant equal to $1$.

	For the moment we are considering lattices in ${\bf C}^n$ of
the form $A(({\bf Z}[i])^n)$, $A \in GL({\bf C}^n)$.  It is natural to
look at these lattices up to unitary equivalence, which is to say that
two lattices $L_1$, $L_2$ are equivalent if there is a unitary linear
transformation $T$ on ${\bf C}^n$ such that $T(L_1) = L_2$.  This
leads to an equivalence relation on $GL({\bf C}^n)$, in which two
invertible linear transformations $A_1$, $A_2$ on ${\bf C}^n$ are
considered to be equivalent if there is a unitary linear
transformation $T$ on ${\bf C}^n$ such that $A_2 = T \, A_1$.

	The quotient of $GL({\bf C}^n)$ by this equivalence relation
is denoted
\begin{equation}
	U({\bf C}^n) \backslash GL({\bf C}^n).
\end{equation}
This quotient space can be identified with the space of self-adjoint
linear transformations on ${\bf C}^n$ which are positive definite,
through the mapping
\begin{equation}
	A \in GL({\bf C}^n) \mapsto A^* \, A.
\end{equation}
That is, for each element $A$ of $GL({\bf C}^n)$, the product
$A^* \, A$ is a self-adjoint linear transformation on ${\bf C}^n$
which is positive-definite,
\begin{equation}
	A_1^* \, A_1 = A_2^* \, A_2
\end{equation}
for two elements $A_1$, $A_2$ of $GL({\bf C}^n)$ if and only if
$A_2 = T \, A_1$ for some unitary linear transformation $T$
on ${\bf C}^n$, and every self-adjoint linear transformation
on ${\bf C}^n$ can be expressed as $A^* \, A$ for some invertible
linear transformation $A$ on ${\bf C}^n$, and in fact has a unique
self-adjoint positive-definite square root.

	Similarly, one can consider two elements $B_1$, $B_2$ of
$SL({\bf C}^n)$ to be equivalent when there is a linear transformation
$U$ in the special unitary group $SU({\bf C}^n)$ such that $A_2 = U \,
A_1$.  The quotient $SU({\bf C}^n) \backslash SL({\bf C}^n)$ can be
identified with the space of self-adjoint linear transformations on
${\bf C}^n$ which are positive-definite and have determinant $1$,
through the same mapping as before.  Also, it will be convenient
to restrict our attention for the moment to lattices $L$ of the form
$B(({\bf Z}[i])^n)$ for some $B \in SL({\bf C}^n)$, which is just
a modest additional normalization.

	Let us write $\Sigma({\bf C}^n)$ for the subgroup of $SL({\bf
C}^n)$ of linear transformations whose associated $n \times n$
matrices, with respect to the standard basis for ${\bf C}^n$, have
integer entries, which implies that the matrices associated to their
inverses also have integer entries.  Thus $B(({\bf Z}[i])^n) = ({\bf
Z}[i])^n$ when $B \in \Sigma({\bf R}^n)$, and conversely $B \in
SL({\bf C}^n)$ and $B(({\bf Z}[i])^n) = ({\bf Z}[i])^n$ implies that
$B \in \Sigma({\bf C}^n)$.  The quotient $SL({\bf C}^n) / \Sigma({\bf
C}^n)$ represents the space of lattices that we are considering at the
moment, the double quotient $SU({\bf C}^n) \backslash SL({\bf C}^n) /
\Sigma({\bf C}^n)$ represents the space of these lattices modulo
equivalence under special unitary transformations, and this double
quotient can also be identified with the quotient of the space of
self-adjoint positive-definite linear transformations on ${\bf C}^n$
with determinant $1$ by the action of $\Sigma({\bf C}^n)$ defined by
$P \mapsto B^* \, P \, B$, $B \in \Sigma({\bf C}^n)$.

	Now let us look at general lattices in ${\bf C}^n$, under the
equivalence relation in which two lattices $L_1$, $L_2$ are considered
to be equivalent if there is an invertible complex-linear
transformation $A$ on ${\bf C}^n$ such that $A(L_1) = L_2$.  This
leads to an equivalence relation on $\mathcal{L}^*({\bf R}^{2n}, {\bf
C}^n)$, in which two elements of $\mathcal{L}^*({\bf R}^{2n}, {\bf
C}^n)$ are considered to be equivalent if one can be written as the
composition of an invertible complex-linear transformation on ${\bf
C}^n$ with the other element of $\mathcal{L}^*({\bf R}^{2n}, {\bf
C}^n)$.  In other words, we look at the action of $GL({\bf C}^n)$ on
$\mathcal{L}^*({\bf R}^{2n}, {\bf C}^n)$ by post-composition.

	Actually, it is more convenient to consider
$\mathcal{L}_1^*({\bf R}^{2n}, {\bf C}^n)$, which we define to be the
subset of $\mathcal{L}^*({\bf R}^{2n}, {\bf C}^n)$ consisting of
invertible real-linear transformations from ${\bf R}^{2n}$ to ${\bf
C}^n$ such that the image of the first $n$ standard basis vectors in
${\bf R}^{2n}$ are linearly-independent over the complex numbers as
$n$ vectors in ${\bf C}^n$.  This restriction is not too serious, and
indeed we can describe the lattices in ${\bf C}^n$ as images of ${\bf
Z}^{2n}$ under mappings in $\mathcal{L}_1^*({\bf R}^{2n}, {\bf C}^n)$.
In other words, if we start with a lattice $L$ given as the image of
${\bf Z}^{2n}$ under an element of $\mathcal{L}^*({\bf R}^n, {\bf
C}^n)$, we can rewrite it as the image of ${\bf Z}^{2n}$ under a
linear transformation in $\mathcal{L}_1^*({\bf R}^{2n}, {\bf C}^n)$ by
pre-composing the initial linear transformation from ${\bf R}^{2n}$ to
${\bf C}^n$ with an invertible linear transformation on ${\bf R}^{2n}$
which permutes the standard basis vectors in a suitable way.

	To deal with the action of $GL({\bf C}^n)$ by
post-composition, we can restrict ourselves to $\mathcal{L}^{**}({\bf
R}^{2n}, {\bf C}^n)$, which we define to be the space of invertible
real-linear transformations from ${\bf R}^{2n}$ to ${\bf C}^n$ such
that the images of the first $n$ standard basis vectors in ${\bf
R}^{2n}$ are the $n$ standard basis vectors in ${\bf C}^n$, and in the
same order.  In other words, if we identify ${\bf R}^{2n}$ with the
Cartesian product ${\bf R}^n \times {\bf R}^n$, then these are the
invertible real-linear transformations from ${\bf R}^n \times {\bf
R}^n$ onto ${\bf C}^n$ with the property that on ${\bf R}^n \times
\{0\}$ they coincide with the standard embedding of ${\bf R}^n$ into
${\bf C}^n$.  This exactly compensates for the action of $GL({\bf C}^n)$
by post-composition, since for any collection $v_1, \ldots, v_n$
of linearly-independent vectors in ${\bf C}^n$ there is a unique
$A \in GL({\bf C}^n)$ such that $A(v_1), \ldots, A(v_n)$
are the standard basis vectors in ${\bf C}^n$, in order.

	We can identify $\mathcal{L}^{**}({\bf R}^{2n}, {\bf C}^n)$
with an open subset of $\mathcal{L}({\bf R}^n, {\bf C}^n)$.  That is,
elements of $\mathcal{L}^{**}({\bf R}^{2n}, {\bf C}^n)$ can be
identified with linear transformations from ${\bf R}^n \times {\bf
R}^n$ into ${\bf C}^n$, and these linear transformations are
determined by what they do on $\{0\} \times {\bf R}^n$, since their
behavior on ${\bf R}^n \times \{0\}$ is fixed by definition.  We can
think of elements of $\mathcal{L}({\bf R}^n, {\bf C}^n)$ as being
written as $A + i \, B$, where $A$, $B$ are linear transformations on
${\bf R}^n$, and one can check that the elements of
$\mathcal{L}^{**}({\bf R}^{2n}, {\bf C}^n)$ correspond exactly to
elements of $\mathcal{L}({\bf R}^n, {\bf C}^n)$ of the form $A + i \,
B$, where $A$, $B$ are linear transformations on ${\bf R}^n$ and $B$
is invertible.

	To be more precise, it is helpful to think in terms of 
real-linear mappings on ${\bf C}^n$, which can be written as
\begin{equation}
	T(x + i \, y) = E_1(x) + E_2(y) + i (E_3(x) + E_4(y)),
\end{equation}
where $x, y \in {\bf R}^n$.  The passage to $\mathcal{L}_1^*({\bf
R}^{2n}, {\bf C}^n)$ can be expressed in these terms as the
restriction to invertible real-linear transformations $T$ on ${\bf
C}^n$ of the form
\begin{equation}
	T(x + i \, y) = x + A(y) + i \, B(y),
\end{equation}
where $A$, $B$ are linear transformations on ${\bf R}^n$.  The
condition of invertibility of $T$ is equivalent to the invertibility
of $B$ on ${\bf C}^n$.

	Another way to look at real-linear mappings on ${\bf C}^n$
is as mappings of the form
\begin{equation}
	T(z) = M(z) + \overline{N(z)},
\end{equation}
where $z \in {\bf C}^n$, $M$ and $N$ are complex-linear mappings on
${\bf C}^n$, and for $w \in {\bf C}^n$, $\overline{w}$ is the element
of ${\bf C}^n$ whose coordinates are the complex-conjugates of the
coordinates of $w$.

	Invertibility of $T$ is a bit tricky, and as an
important special case, it is natural to restrict our attention to
mappings $T$ as above for which $M$ majorizes $N$ in the sense that
\begin{equation}
	|N(z)| < |M(z)|
\end{equation}
for $z \in {\bf C}^n$, $z \ne 0$, where $|w|$ denotes the standard
Euclidean norm of $w \in {\bf C}^n$.  To factor out the action of
$GL({\bf C}^n)$ by post-composition, we can restrict our attention to
real-linear transformations $T$ of the form
\begin{equation}
	T(z) = z + \overline{E(z)},
\end{equation}
where $E$ is a complex-linear transformation on ${\bf C}^n$ with
operator norm strictly less than $1$, which is equivalent to saying
that $E^* \, E < I$.  This has nice features when we think of the
image of the standard integer lattice $({\bf Z}[i])^n$ under $T$, with
points in the image being reasonably-close to their counterparts in
the original lattice.

	The $n = 1$ case is quite instructive.  We can write a
real-linear transformation $T$ on ${\bf C}$ as
\begin{equation}
	T(x + i \, y) = a \, x + i \, b \, y
\end{equation}
for $x, y \in {\bf R}$, where $a$, $b$ are complex numbers, and when
$T$ is invertible we can rewrite this as
\begin{equation}
	T(x + i \, y) = a (x + i \, c \, y),
\end{equation}
where $a$, $c$ are complex numbers with $a \ne 0$ and $c$ having
nonzero imaginary part.  Alternatively, we can write a real-linear
transformation $T$ on ${\bf C}$ as $T(z) = \alpha \, z + \beta \,
\overline{z}$ with $\alpha, \beta \in {\bf C}$, and where $T$ is
invertible if and only if $|\alpha| \ne |\beta|$, and when
$|\alpha| > |\beta|$ this can be rewritten as
\begin{equation}
	T(z) = \theta (z + \mu \, \overline{z}),
\end{equation}
where $\theta$ is a nonzero complex number and $\mu$ is a complex
number such that $|\mu| < 1$.

\end{document}